\documentclass[a4paper,11pt]{amsart}
\usepackage{amssymb,amscd}
\usepackage[cp1251]{inputenc}
\usepackage[T2A]{fontenc}
\theoremstyle{plain}

\newtheorem{theo}{Theorem}[section]
\newtheorem{prop}[theo]{Proposition}
\newtheorem{lemm}[theo]{Lemma}
\newtheorem{coro}[theo]{Corollary}
\newtheorem{prob}[theo]{Problem}

\theoremstyle{definition}
\newtheorem*{defi}{Definition}
\newtheorem{exam}[theo]{Example}

\theoremstyle{remark}
\newtheorem*{rema}{Remark}

\numberwithin{equation}{section}

\newcommand{\field}[1]{\mathbb{#1}}
\def\C{\field{C}}
\newcommand{\R}{\field{R}}
\newcommand{\Q}{\field{Q}}
\newcommand{\Z}{\field{Z}}

\def\N{\mathit N}
\def\T{\mathit T}

\def\L{\varLambda}

\def\BQ{\mbox{\textit{BQ}}}

\def\f{\varphi}

\DeclareMathOperator{\lk}{lk}
\DeclareMathOperator{\st}{st}
\DeclareMathOperator{\chf}{cf}

\def\ge{\geqslant}

\def\le{\leqslant}

\begin{document}

\title{Semifree circle actions, Bott towers, and quasitoric manifolds}

\author{Mikiya Masuda}
\author{Taras Panov}\thanks{The second author gratefully
acknowledges the support from the Japanese Society for the
Promotion of Science (grant~P05296), the Russian Foundation for
Basic Research (grant~05-01-01032) and P.~Deligne's 2004 Balzan
prize in mathematics.}
\address{Department of Mathematics, Osaka City University,
Sumiyoshi-ku, Osaka 558-8585, Japan}
\address{Department of Geometry and Topology, Faculty of Mathematics and Mechanics, Moscow State University,
Leninskiye Gory, Moscow 119992, Russia}
\email{masuda@sci.osaka-cu.ac.jp} \email{tpanov@mech.math.msu.su}


\begin{abstract}
A \emph{Bott tower} is the total space of a tower of fibre bundles
with base $\mathbb C P^1$ and fibres~$\mathbb C P^1$. Every Bott
tower of height $n$ is a smooth projective toric variety whose
moment polytope is combinatorially equivalent to an $n$-cube. A
circle action is \emph{semifree} if it is free on the complement
to fixed points. We show that a (quasi)toric manifold (in the
sense of Davis--Januszkiewicz) over an $n$-cube with a semifree
circle action and isolated fixed points is a Bott tower. Then we
show that every Bott tower obtained in this way is topologically
trivial, that is, homeomorphic to a product of 2-spheres. This
extends a recent result of Ilinskii, who showed that a smooth
compact toric variety with a semifree circle action and isolated
fixed points is homeomorphic to a product of 2-spheres, and makes
a further step towards our understanding of a problem motivated by
Hattori's work on semifree circle actions. Finally, we show that
if the cohomology ring of a quasitoric manifold (or Bott tower) is
isomorphic to that of a product of 2-spheres, then the manifold is
homeomorphic to the product.
\end{abstract}

\maketitle

\section{Introduction}
In their study of symmetric spaces Bott and
Samelson~\cite{bo-sa58} introduced a family of complex manifolds
obtained as the total spaces of iterated bundles over $\mathbb C
P^1$ with fibre $\mathbb C P^1$. Grossberg and
Karshon~\cite{gr-ka94} showed these manifolds to carry an
algebraic torus action, therefore constituting an important family
of smooth projective toric varieties, and called them \emph{Bott
towers}. The study of Bott towers has been significantly developed
by Civan and Ray in~\cite{ci-ra05} by enumerating the invariant
stably complex structures and calculating their complex and real
$K$-theory rings, and cobordisms.

An action of a group is called \emph{semifree} if it is free in
the complement to fixed points. A particularly interesting class
of \emph{Hamiltonian} semifree circle actions was studied by
Hattori, who proved in~\cite{hatt92} that a compact symplectic
manifold $M$ with a semifree Hamiltonian $\mathbb S^1$-action with
non-empty isolated fixed point set has the same cohomology ring
and the same Chern classes as $S^2\times\ldots\times S^2$, thus
imposing a severe topological restriction on the topological
structure of the manifold. Hattori's results were further extended
by Tolman and Weitsman, who showed in~\cite{to-we00} that a
semifree symplectic $\mathbb S^1$-action with non-empty isolated
fixed point set is automatically Hamiltonian, and the
\emph{equivariant} cohomology ring and Chern classes of $M$ also
agree with those of $S^2\times\ldots\times S^2$. In dimensions up
to 6 it is known that a symplectic manifold with an $\mathbb
S^1$-action satisfying the properties above is homeomorphic to a
product of 2-spheres, but in higher dimensions the topological
classification remains open.

Ilinskii considered in~\cite{ilin06} an algebraic version of the
above question on semifree symplectic $\mathbb S^1$-actions.
Namely, he conjectured that a smooth compact complex algebraic
variety $X$ with a semifree action of the algebraic 1-torus $\C^*$
with positive number of isolated fixed points is homeomorphic to
$\C P^1\times\ldots\times \C P^1$. The algebraic and symplectic
versions of the conjecture are related via the common subclass of
projective varieties; a smooth projective variety is a symplectic
manifold. Ilinskii proved the \emph{toric} version of his
algebraic conjecture, namely, the case when $X$ is a \emph{toric
manifold} (a non-singular compact \emph{toric variety}) and the
semifree 1-torus is a subgroup of the acting torus (of dimension
$\dim_\C X$). The first step of Ilinskii's argument was to show
that if $X$ admits a semifree subcircle with isolated fixed
points, then the corresponding fan is combinatorially equivalent
to the fan over the faces of a crosspolytope. A result of
Dobrinskaya~\cite{dobr01} implies that such $X$ is a Bott tower,
and this was the starting point in our study of semifree circle
actions on Bott towers and related classes of manifolds with torus
action, such as the \emph{quasitoric manifolds}. This class of
manifolds was introduced by Davis and Januszkiewicz
in~\cite{da-ja91} (they used the name ``toric manifolds'', which
we prefer to save for ``non-singular compact toric varieties'' as
above). A quasitoric manifold is a compact manifold $M$ of
dimension $2n$ with a locally standard action of an
$n$-dimensional torus $\mathbb T^n=\mathbb S^1\times\ldots\times
\mathbb S^1$ whose orbit space is a simple polytope~$P$. Recently
quasitoric manifolds have attracted substantial interest in the
emerging field of ``toric topology''; we review their construction
in Section~\ref{quasi}, for a more detailed account
see~\cite[Ch.~5]{bu-pa02}.

A projective toric manifold with the moment polytope $P$ is a
quasitoric manifold over $P$, and a Bott tower is a toric manifold
with the moment polytope combinatorially equivalent to a cube (or
simply a \emph{toric manifold over a cube}). We therefore have the
following hierarchy of classes of manifolds $M$ with an action of
$\mathbb T^n$:
\begin{multline}\label{hiera}
  \text{Bott towers }
  \subset\text{ Toric manifolds over cubes}\\
  \subset\text{ Quasitoric manifolds over cubes.}
\end{multline}
By the above mentioned result of Dobrinskaya~\cite{dobr01}, the
first inclusion above is in fact an identity (we explain this in
Corollary~\ref{torbt}). We proceed in Sections~\ref{semif}
and~\ref{tocla} by obtaining two results relating semifree circle
actions on Bott towers, their topological structure, and
cohomology rings. In Theorem~\ref{btpro} we show that if a Bott
tower admits a semifree $\mathbb S^1$-action with isolated fixed
points, then it is $\mathbb S^1$-equivariantly homeomorphic to a
product of 2-spheres. We also show in Theorem~\ref{trivc} that a
Bott tower with the cohomology ring isomorphic to that of a
product of spheres is homeomorphic to the product. Both results
are then extended to a much more general class of quasitoric
manifolds over cubes (Theorems~\ref{qtmsf} and~\ref{qtcpn}
respectively), which also allows us to deduce Ilinskii's result on
semifree 1-dimensional algebraic torus actions on toric varieties
(Corollary~\ref{ilins}).

Since a cohomology isomorphism implies a homeomorphism in the case
of product of spheres,
we may ask a question of whether the cohomology ring detects the
homeomorphism type of a Bott tower or quasitoric manifold in
general. Some progress in this direction has been achieved in
upcoming work~\cite{c-m-s??} in the case of quasitoric manifolds
over a product of simplices, which can be considered as an
intermediate stage between quasitorics over a cube and the full
generality.

It would be interesting to obtain smooth analogues of our
classification results. In the case of Bott towers our
homeomorphisms are actually diffeomorphisms (see
Theorem~\ref{btpro} or Theorem~\ref{trivc}), but some of our key
arguments with quasitoric manifolds do not work in the smooth
category. Although quasitoric manifolds are necessarily
smooth~\cite[p.~421]{da-ja91}, the main problem here is that the
original Davis and Januszkiewicz's classification
result~\cite[Prop.~1.8]{da-ja91} establishes only an equivariant
\emph{homeo}morphism between a quasitoric manifold and the
canonical model determined by the polytope and the characteristic
function. As the consequence, we don't know if there are exotic
equivariant smooth structures, even on 4-dimensional quasitoric
manifolds. (A canonical equivariant smooth structure, coinciding
with the standard one in the toric case, is described
in~\cite[\S4]{b-p-r07}.)

All the cohomology groups in the paper are taken with $\Z$
coefficients, unless otherwise specified.

The authors are grateful to Natalia Dobrinskaya and Dong Youp Suh
for informal discussions of quasitoric manifolds over cubes and
producs of simplices, and apologise for the improper references
to~\cite{dobr01} in the first version of the text. We also thank
Takahiko Yoshida, who pointed out the difficulty in the smooth
category for quasitoric manifolds. The second author wishes to
thank Nigel Ray for introducing him to the study of Bott towers
and illuminating discussions of the subject. Second author's
thanks also go to Dmitri Timashev for drawing his attention to
work of Ilinskii~\cite{ilin06}, which encouraged us to consider
several related problems on semifree circle actions.

\section{Bott towers}
We briefly review the definitions here, following~\cite{gr-ka94}
and~\cite{ci-ra05}, where a reader may find a much more detailed
account of the history and applications of Bott towers.

\begin{defi}
A \emph{Bott tower} of height $n$ is a sequence of manifolds
$(B^{2k}\colon k\le n)$ such that $B^2=\C P^1$ and
$B^{2k}=P(\underline{\C}\oplus\xi_{k-1})$ for $1<k\le n$ where
$P(\cdot)$ denotes complex projectivisation, $\xi_{k-1}$ is a
complex line bundle over~$B^{2(k-1)}$ and $\underline{\C}$ is a
trivial line bundle. In particular, we have a bundle $B^{2k}\to
B^{2(k-1)}$ with fibre~$\C P^1$.
\end{defi}

We shall also use the same name ``Bott tower'' for the last stage
$B^{2n}$ in the sequence; it follows from the definition that
$B^{2n}$ is a complex manifold obtained as the total space of an
iterated bundle with fibre~$\C P^1$. Bott towers of (real)
dimension 4 are known as \emph{Hirzebruch surfaces}.

The standard results on the cohomology of projectivised bundles
lead to the following description of the cohomology ring of a Bott
tower (see e.g.~\cite[Cor.~2.9]{ci-ra05}).

\begin{lemm}\label{cohbt}
The cohomology of $B^{2k}$ is a free module over $H^*(B^{2(k-1)})$
on generators 1 and $u_k$, which have dimensions 0 and 2
respectively; the ring structure is determined by the single
relation
\[
  u_k^2=c_1(\xi_{k-1})u_k,
\]
and the restriction of $u_k$ to the fibre~$\C P^1\subset B^{2k}$
is the first Chern class of the canonical line bundle over $\C
P^1$.
\end{lemm}

We set $u_1$ to be the canonical generator of $H^2(\C P^1)$ (the
first Chern class of the canonical line bundle). The Bott tower is
therefore determined by the list of integers $(a_{ij}\colon 1\le
i<j\le n)$, where
\begin{equation}\label{crobt}
  u_k^2=\sum_{i=1}^{k-1}a_{ik}u_iu_k, \qquad 1\le k\le n.
\end{equation}
The cohomology ring of $B^{2n}$ is the quotient of
$\Z[u_1,\ldots,u_n]$ by relations \eqref{crobt}.

It is convenient to organise the integers $a_{ij}$ into an
$n\times n$ upper triangular matrix,
\begin{equation}\label{amatr}
  A=\begin{pmatrix}
    -1 & a_{12} & \cdots & a_{1n}\\
    0  & -1     & \cdots & a_{2n}\\
    \vdots & \vdots & \ddots & \vdots\\
    0 & 0 & \cdots & -1
  \end{pmatrix}.
\end{equation}

\begin{exam}\label{2dbto}
When $n=2$, the Bott tower $B^4$ is determined by a single line
bundle $\xi_1$ over $\C P^1$. We have $\xi_1=\gamma^m$ for some
$m\in\Z$ where $\gamma$ is the canonical line bundle over $B^2=\C
P^1$, and so the cohomology ring is determined by the relations
$u_1^2=0$ and $u_2^2=mu_1u_2$. It is well-known that
\[
  P(\underline{\C}\oplus\gamma^m)\cong
  P(\underline{\C}\oplus\gamma^{m'})
  \text{ if and only if }m\equiv m'
  \pmod 2,
\]
where $\cong$ means ``diffeomorphic''. The proof goes as follows.
We note that $P(E)\cong P(E\otimes\eta)$ for any complex line
bundle $\eta$. Suppose $m\equiv m' \pmod 2$. Then $m'-m=2\ell$ for
some $\ell\in \Z$ and we have diffeomorphisms
\[
P(\underline{\C}\oplus\gamma^m)\cong
P((\underline{\C}\oplus\gamma^m)\otimes\gamma^{\ell})=
P(\gamma^\ell\oplus\gamma^{m+\ell}).
\]
Here both $\gamma^\ell\oplus\gamma^{m+\ell}$ and
$\underline{\C}\oplus\gamma^{m'}$ are over $\C P^1$ and have equal
first Chern class, so they are isomorphic. Hence the last space
above is $P(\underline{\C}\oplus\gamma^{m'})$. On the other hand,
it is not difficult to see that if
$H^*(P(\underline{\C}\oplus\gamma^m)) \cong
H^*(P(\underline{\C}\oplus\gamma^{m'}))$ as rings, then $m\equiv
m'\pmod 2$.
\end{exam}

This example shows that the cohomology ring determines the
topological type of a Bott tower $B^{2n}$ for $n=2$. A
case-to-case analysis based on a classification result
of~\cite[\S3]{dobr01} shows that this is also the case for $n=3$.
So we may ask the following question.

\begin{prob}\label{cohcl}
Are Bott towers $B_1^{2n}$ and $B_2^{2n}$ homeomorphic if
$H^*(B_1^{2n})\cong H^*(B_2^{2n})$ as rings?
\end{prob}

We investigate this question further in Section~\ref{tocla}, where
a partial answer is given.

\section{Quasitoric manifolds}\label{quasi}
Davis and Januszkiewicz introduced in~\cite{da-ja91} a class of
$2n$-dimensional manifolds $M$ with an action of $n$-dimensional
torus $T$. They required the action to be \emph{locally standard}
(locally isomorphic to the standard $T$-representation in $\mathbb
C^n$) and the quotient space $M/T$ to be homeomorphic to a
\emph{simple} $n$-dimensional polytope $P$, so there is a
projection $\pi\colon M\to P$ whose fibres are orbits of the
action. Davis and Januszkiewicz describe their manifolds as {\it
toric}; more recently, the term {\it quasitoric} has been adopted,
to avoid confusion with the non-singular compact toric varieties
of algebraic geometry. We follow this convention below, and refer
to such $M$ as {\it quasitoric manifolds}, saving the term ``toric
manifolds'' to describe non-singular compact toric varieties. We
note that a projective toric manifold is a quasitoric manifold;
for more discussion of the relationship between the two classes
see~\cite[Ch.~5]{bu-pa02}.

Let $m$ denote the number of \emph{facets} (codimension-one faces)
of $P$; we order the facets so that the first $n$ of them meet at
a vertex. We denote the facets by $F_i$ for $1\le i\le m$, and
denote $\mathcal F$ the set of all facets. The preimage
$\pi^{-1}(F_i)$ is a connected codimension-two submanifold of $M$,
fixed pointwise by a circle subgroup of $\T$. We denote it by
$M_i$ and refer to it as the \emph{characteristic submanifold}
corresponding to $F_i$, for $1\le i\le m$. An
\emph{omniorientation}~\cite{bu-ra01} of $M$ consists of a choice
of orientation for $M$ and for each characteristic submanifold.

Let $\N$ denote the integer lattice of one-parameter circle
subgroups in $\T$, so we have $\N\cong\Z^n$. Given a
characteristic submanifold $M_j$, we denote by $\lambda_j$ a
primitive vector in $\N$ that spans the circle subgroup
$T_{M_j}\subset\T$ fixing $M_j$. The vector $\lambda_j$ is
determined up to sign. The correspondence $\lambda\colon
F_j\mapsto\lambda_j$ was called in~\cite{da-ja91} the
\emph{characteristic function} corresponding to $M$.

The omniorientation allows us to interpret the characteristic
function as a linear map $\lambda\colon\Z^{\mathcal F}\to\N$. This
is done as follows. First, we notice that an action of a
one-parameter circle subgroup of $T$ determines an orientation for
the normal bundle $\nu_j$ of the embedding $M_j\subset M$. The
omniorientation of $M$ also provides an orientation for $\nu_j$ by
means of the following decomposition of the tangent bundle:
\[
  \tau(M)|_{M_j}=\tau(M_j)\oplus\nu_j.
\]
Now we choose the primitive vectors $\lambda_j$ so that the two
orientations of $\nu_j$ coincide for $1\le j\le m$.

In general, there is no canonical choice of an omniorientation for
$M$. However, if $M$ admits a \emph{$T$-equivariant almost complex
structure}, then a choice of such a structure provides a canonical
way of orienting $M$ and normal bundles $\nu_j$ for $1\le j\le m$,
thereby specifying an omniorientation \emph{associated} with the
equivariant almost complex structure. In what follows we shall
always choose the associated omniorientation if $M$ is
equivariantly almost complex (in particular, if $M$ is a toric
manifold); otherwise we shall fix an arbitrary omniorientation.

By definition, the characteristic function satisfies the
\emph{non-singularity condition}:
$\lambda_{j_1},\ldots,\lambda_{j_n}$ is a basis of $\N$ whenever
the intersection $F_{j_1}\cap\ldots\cap F_{j_n}$ is non-empty. So
we may use the vectors $\lambda_1,\ldots,\lambda_n$ to define a
basis for $\N$, thereby identifying it with $\Z^n$, and represent
the map $\lambda$ by an $n\times m$ integral matrix of the form
\begin{equation}\label{lamat}
\L\;=\;\begin{pmatrix}
  1&0&\ldots&0&\lambda_{1,n+1}&\ldots&\lambda_{1,m}\\
  0&1&\ldots&0&\lambda_{2,n+1}&\ldots&\lambda_{2,m}\\
  \vdots&\vdots&\ddots&\vdots&\vdots&\ddots&\vdots\\
  0&0&\ldots&1&\lambda_{n,n+1}&\ldots&\lambda_{n,m}
\end{pmatrix}.
\end{equation}
It is often convenient to partition $\L$ as
$\left(E\;|\;\L_\star\right)$, where $E$ is a unit matrix and
$\L_\star$ is $n\times(m-n)$. For any other vertex
$F_{j_1}\cap\ldots\cap F_{j_n}$ the corresponding columns
$\lambda_{j_1},\ldots,\lambda_{j_n}$ form a basis for $\Z^n$, and
have determinant $\pm1$. We refer to~\eqref{lamat} as the
\emph{refined form} of characteristic matrix $\L$, and call
$\L_\star$ its \emph{reduced submatrix}.

Having chosen a basis for $\N$, we may identify our torus $\T$
with the standard product of unit circles in $\C^n$:
\begin{equation}\label{torus}
  \mathbb T^n=\{(e^{2\pi i\f_1},\ldots,e^{2\pi i\f_n})\in\C^n\},
\end{equation}
where $(\f_1,\ldots,\f_n)$ runs over $\R^n$. We shall also denote
a generic point of $\mathbb T^n$ by $(t_1,\ldots,t_n)$. The circle
subgroup fixing $M_j$ can now be written as
\begin{equation}\label{tmj}
  T_{M_j}=\{(t^{\lambda_{1j}},\ldots,t^{\lambda_{nj}})=
  (e^{2\pi i\lambda_{1j}\f},\ldots,e^{2\pi
  i\lambda_{nj}\f})\in \mathbb T^n\},\quad 1\le j\le m
\end{equation}
where $\f\in\R$ and $t=e^{2\pi i\f}$.

\begin{rema}
Not every toric manifold $X$ is a quasitoric manifold, as the
quotient $X/\T$ may fail to be a simple polytope (although it is a
simple polytope when $X$ is projective). Nevertheless, $X$ has
characteristic submanifolds $X_j$ ($T$-invariant divisors), and
there is a canonical omniorientation induced from the complex
structures on $X$ and $X_j$. Therefore, the characteristic matrix
$\L$ and its reduced submatrix $\L_\star$ are defined for every
toric manifold $X$. The vectors $\lambda_j$ are the primitive
vectors along the edges of the fan corresponding to~$X$.
\end{rema}

Let $v_j$ denote the class in $H^2(M)$ dual to the fundamental
class of $M_j$, $1\le j\le m$. According to
\cite[Th.~4.14]{da-ja91}, the ring $H^*(M)$ is generated by
$v_1,\ldots,v_m$, modulo two sets of relations. The first set is
formed by the monomial relations which arise from the
\emph{Stanley--Reisner ideal} of $P$; the second set consists of
linear relations determined by the characteristic matrix:
\begin{equation}\label{redvi}
  v_i\;=\;-\lambda_{i,n+1}v_{n+1}-\ldots-\lambda_{i,m}v_m
  \text{ \ for \ }1\le i\le n.
\end{equation}
It follows that $v_{n+1},\ldots,v_m$ suffice to generate $H^*(M)$
multiplicatively.

Two quasitoric manifolds $M_1$ and $M_2$ are said to be
\emph{weakly $\T$-equivariantly homeomorphic} (or simply
\emph{weakly $\T$-homeomorphic}) if there is an automorphism
$\theta\colon\T\to\T$ and a homeomorphism $f\colon M_1\to M_2$
such that $f(t\cdot x)=\theta(t)\cdot f(x)$ for every $t\in\T$ and
$x\in M_1$. If $\theta$ is identity, then $M_1$ and $M_2$ are
\emph{$\T$-homeomorphic}. Following Davis and Januszkiewicz, we
say that two quasitoric manifolds $M_1$ and $M_2$ over the same
$P$ are \emph{equivalent} if there is a weak $T$-homeomorphism
$f\colon M_1\to M_2$ covering the identity on $P$.
By~\cite[Prop.~1.8]{da-ja91}, a quasitoric manifold $M$ over $P$
is determined up to equivalence by its characteristic
function~$\lambda$. This follows from the ``basic construction''
providing a quasitoric manifold $M(\lambda)$, which depends only
on $P$ and $\lambda$, together with a $\T$-equivariant
homeomorphism $M(\L)\to M$ covering the identity on~$P$.

Let $\chf(P)$ denote the set of characteristic functions on the
facets of $P$, that is, the set of maps $\lambda\colon\mathcal
F\to N$ satisfying the non-singularity condition. The group
$GL(\N)\cong GL(n,\Z)$ of automorphisms of the lattice $\N$ acts
on the set $\chf(P)$ from the left (an automorphism
$g\colon\N\to\N$ acts by composition $\lambda\mapsto
g\cdot\lambda$). Since automorphisms of the lattice $\N$
correspond to automorphisms of the torus $\T$, there is a
one-to-one correspondence
\begin{equation}\label{lcset}
  GL(\N)\backslash\chf(P)\longleftrightarrow
  \{\text{equivalence classes of $M$ over $P$}\}.
\end{equation}
One may assign an $n\times m$-matrix $\Lambda$ to an element
$\lambda\in\chf(P)$ by ordering the facets and choosing a basis
for $\N$, as we did above. A choice of matrix $\L$ in the refined
form~\eqref{lamat} can now be regarded as a choice of a specific
representative of the left coset class in
$GL(\N)\backslash\chf(P)$. If a characteristic matrix is given in
an unrefined form $\L=(A\;|\;B)$, where $A$ is $n\times n$ and $B$
is $n\times(m-n)$, then the refined representative in its coset
class is given by $(E\;|\;A^{-1}B)$.

A quasitoric manifold over $P$ can be obtained as the quotient of
the \emph{moment-angle manifold} $\mathcal Z_P$ by a freely acting
$(m-n)$-dimensional torus subgroup of $\mathbb T^m$ determined by
the kernel of the characteristic map $\lambda\colon\Z^m\to\N$, see
\cite{da-ja91}, \cite[Ch.~6]{bu-pa02}. The moment-angle manifold
$\mathcal Z_P$ embeds into $\C^m$ as a complete intersection of
$m-n$ real quadratic hypersurfaces~\cite[\S3]{b-p-r07}. It follows
that both $\mathcal Z_P$ and $M$ are necessarily smooth, but we
don't know if this equivariant smooth structure on $M$ is unique
(see the discussion in the Introduction).

We are particularly interested in the case when the quotient
polytope $P=M\!/\,\mathbb T^n$ is the $n$-cube $\mathbb I^n$. Then
$m=2n$ and we shall additionally assume that the facets $F_j$ and
$F_{n+j}$ are opposite (i.e., do not intersect) for $1\le j\le n$.
In the case $P=\mathbb I^n$ the moment-angle manifold is the
product of $n$ three-dimensional spheres, embedded in $\C^{2n}$ as
\[
  \{(z_1,\ldots,z_{2n})\in\C^{2n}\colon|z_j|^2+|z_{n+j}|^2=1\text{ for }1\le
  j\le n\}.
\]
The quotient $(S^3)^n\!/\,\mathbb T^{2n}$ by the coordinatewise
action is a cube $\mathbb I^n$. The $n$-dimensional subtorus
$T(\L)\subset\mathbb T^{2n}$ determined by the kernel of
characteristic map~\eqref{lamat} is given by
\begin{multline}\label{subto}
  (t_1,\ldots,t_n)\mapsto
  (t_1^{-\lambda_{1,n+1}}t_2^{-\lambda_{1,n+2}}\cdots
  t_{n}^{-\lambda_{1,2n}},\ldots,\\
  t_1^{-\lambda_{n,n+1}}t_2^{-\lambda_{n,n+2}}\cdots
  t_{n}^{-\lambda_{n,2n}},
  t_1,t_2,\ldots,t_n).
\end{multline}
It acts freely on $(S^3)^n$, and the quotient $(S^3)^n\!/\,T(\L)$
is the quasitoric manifold $M$ determined by $\L$. The
$n$-dimensional torus $\mathbb T^{2n}\!/\,T(\L)\cong \mathbb T^n$
acts on $(S^3)^n\!/\,T(\L)\cong M$ with quotient $\mathbb I^n$. In
coordinates, $(t_1,\dots,t_n)\in \mathbb T^n$ acts on the
equivalence class $[z_1,\dots,z_{2n}]\in(S^3)^n\!/\,T(\L)$ as
multiplication by $(t_1,\dots,t_n,1,\ldots,1)$.

\begin{prop}\label{almat}
A Bott tower carries a natural torus action turning it into a
quasitoric manifold over a cube with the reduced submatrix
$\L_\star=A^t$, see~\eqref{amatr} and~\eqref{lamat}.
\end{prop}
\begin{proof}
As it is shown in~\cite[Prop.~3.1]{ci-ra05}, the Bott tower
corresponding to~\eqref{amatr} can be obtained as the quotient of
$(S^3)^n$ by an $n$-dimensional subtorus of $\mathbb T^{2n}$
defined by the inclusion
\[
  (t_1,\ldots,t_n)\mapsto
  (t_1,t_1^{-a_{12}}t_2,\ldots,t_1^{-a_{1n}}t_2^{-a_{2n}}\cdots
  t_{n-1}^{-a_{n-1,n}}t_n,t_1,t_2,\ldots,t_n).
\]
It remains to note that this coincides with $T(\L)$
from~\eqref{subto} for $\L_\star=A^t$.
\end{proof}

\begin{rema}
The Stanley--Reisner relations for the $n$-cube are
$v_iv_{i+n}=0$, $1\le i\le n$. These relations together
with~\eqref{redvi} give~\eqref{crobt} for $\L_\star=A^t$ and
$u_i=v_{i+n}$, \ $1\le i\le n$.
\end{rema}

By definition, a Bott tower is a complex manifold. A non-compact
version of~\eqref{subto} was used in~\cite{gr-ka94} to describe
Bott towers as non-singular projective toric varieties. The two
approaches are related in~\cite[\S2]{ci-ra05}.

Given a permutation $\sigma$ of $n$ elements, denote by
$P(\sigma)$ the corresponding $n\times n$ \emph{permutation
matrix}, which has units in positions $(\sigma(i),i)$ for $1\le
i\le n$, and zeros elsewhere. There is an action of the symmetric
group on $n\times n$ matrices by conjugations $A\mapsto
P(\sigma)^{-1}AP(\sigma)$ or, equivalently, by permutations of the
rows and columns of~$A$.

\begin{prop}\label{qtbtc}
A quasitoric manifold $M$ over a cube with reduced submatrix
$\L_\star$ is equivalent to a Bott tower if and only if $\L_\star$
can be conjugated by a permutation matrix to an upper triangular
matrix.
\end{prop}
\begin{proof}
Assume that $\L_\star$ can be conjugated by a permutation matrix
to an upper triangular matrix. Consider the action of $S_n$ on the
set of facets of $\mathbb I^n$ by permuting the pairs of opposite
facets. A reordering of facets corresponds to reordering of
columns in the $n\times2n$ characteristic matrix $\L$, so an
element $\sigma\in S_n$ acts as
\[
  \L\mapsto\L\cdot
  \begin{pmatrix}P(\sigma)&0\\0&P(\sigma)\end{pmatrix}.
\]
This action does not preserve the refined form of $\L$, as
$(E\;|\;\L_\star)$ becomes $(P(\sigma)\;|\;\L_\star P(\sigma))$.
The refined representative in the left coset class~\eqref{lcset}
of the latter matrix is given by $(E\;|\;P(\sigma)^{-1}\L_\star
P(\sigma))$. (In other words, we have to compensate the
permutation of pairs of facets by an automorphism of $\mathbb T^n$
permuting the coordinate subcircles to keep the characteristic
matrix in the refined form.) This implies that the permutation
action on the pairs of opposite facets induces the conjugation
action on the reduced submatrices. So we may assume, up to
equivalence, that $M$ has an upper triangular reduced submatrix
$\L_\star$. The non-singularity condition guarantees that the
diagonal entries of $\L_\star$ are $\pm1$, and we can set all of
them to $-1$ by changing the omniorientation of $M$ if necessary.
Now, $M$ and the Bott tower corresponding to the matrix
$A=\L_\star^t$ have the same characteristic matrices $\L$ by
Proposition~\ref{almat}. Therefore, they are equivalent
by~\cite[Prop.~1.8]{da-ja91}. The opposite statement follows from
Proposition~\ref{almat}.
\end{proof}

It is now clear that not all quasitoric manifolds over a cube are
Bott towers. For example, a 4-dimensional quasitoric manifold over
a square with reduced submatrix
$\L_\star=\begin{pmatrix}-1&-2\\-1&-1\end{pmatrix}$ is not a Bott
tower, as $\L_\star$ cannot be conjugated to an upper triangular
matrix. (The corresponding manifold is homeomorphic to $\C
P^2\mathbin{\#}\C P^2$, and therefore does not even admit a
complex structure~\cite{da-ja91}.)

Given a $k$-element subset $\{i_1,\ldots,i_k\}$ of $n$ elements,
the corresponding \emph{principal minor} of a square $n$-matrix
$A$ is the determinant of the submatrix formed by elements in
columns and rows $i_1,\ldots,i_k$. For Bott towers,
Proposition~\ref{almat} ensures that all principal minors of the
matrix $-\L_\star$ equal 1; while for arbitrary quasitoric
manifolds the non-singularity condition only guarantees that every
principal minor of $\L_\star$ is $\pm1$.

Recall that an upper triangular matrix is \emph{unipotent} if all
its diagonal entries are unit. The following key technical lemma
can be retrieved from the proof of a much more general result of
Dobrinskaya~\cite[Th.~6]{dobr01}. We give a slightly more expanded
proof here for the sake of completeness.

\begin{lemm}[\cite{dobr01}]\label{minor}
Let $R$ be a commutative integral domain with an identity element
$1$, and let $A$ be an $n\times n$ matrix with entries in $R$.
Suppose that every proper principal minor of $A$ is~$1$.  If $\det
A=1$, then $A$ is conjugate by a permutation matrix to a unipotent
upper triangular matrix, and otherwise to a matrix of the form
\begin{equation}\label{exce}
\begin{pmatrix}
1 & b_1 & 0 & \dots &  0\\
0 & 1 & b_2 & \dots &  0\\
\vdots& \vdots  &\ddots &\ddots & \vdots \\
0 & 0 & \dots & 1 &  b_{n-1}\\
b_n & 0 & \dots & 0  & 1
\end{pmatrix}
\end{equation}
where $b_i\not=0$ for every $i$.
\end{lemm}
\begin{proof}
Every diagonal entry of $A$ must be $1$. We say that the $i$-th
row is \emph{elementary} if its $i$-th entry is $1$ and others
are~$0$. Assuming by induction that the theorem holds for matrices
of size $(n-1)$, we deduce that $A$ itself is conjugate to a
unipotent upper triangular matrix if and only if it contains an
elementary row. Denote by $A_i$ the square matrix of size $(n-1)$
obtained by removing the $i$-th column and row from~$A$.

We may assume by induction that $A_n$ is a unipotent upper
triangular matrix. Now we apply the induction assumption to $A_1$.
The permutation of rows and columns turning $A_1$ into a unipotent
upper triangular matrix turns $A$ into an ``almost'' unipotent
upper triangular matrix; the latter may have only one non-zero
entry below the diagonal, and it must be in the first column. If
this non-zero entry is not $a_{n1}$, then the $n$-th row of $A$ is
elementary, and $A$ is conjugate to a unipotent upper triangular
matrix. Therefore, we may assume that
\begin{equation*}
  A=
  \begin{pmatrix}
    1 & * & * & \dots &  *\\
    0 & 1 & * & \dots &  *\\
    \vdots& \vdots &\ddots &\ddots & \vdots \\
    0 & 0 & \dots & 1 &  b_{n-1}\\
    b_n & 0 & \dots & 0  & 1
  \end{pmatrix},
\end{equation*}
where $b_{n-1}\ne0$ and $b_n\ne0$ (otherwise $A$ has an elementary
row). Now let $a_{1j_1}$ be the last non-zero non-diagonal entry
in the first row of~$A$. If $A$ does not have an elementary row,
the we may define inductively $a_{j_ij_{i+1}}$ as the last
non-zero non-diagonal entry in the $j_i$-th row of~$A$. Clearly,
we have
\[
  1<j_1<\ldots<j_i<j_{i+1}<\ldots<j_k=n
\]
for some $k<n$. Now, if $j_i=i+1$ for $i=1,\ldots,n-1$, then $A$
is matrix~\eqref{exce} with $b_i=a_{j_{i-1}j_i}$, \
$i=1,\ldots,n-1$. Otherwise, the submatrix
\begin{equation*}
  S=
  \begin{pmatrix}
    1 & a_{1j_1} & 0 & \dots &  0\\
    0 & 1 & a_{j_1j_2} & \dots &  0\\
    \vdots& \vdots &\ddots &\ddots & \vdots \\
    0 & 0 & \dots & 1 & a_{j_{k-1}n} \\
    b_n & 0 & \dots & 0  & 1
  \end{pmatrix}
\end{equation*}
of $A$, formed by columns and rows $1,j_1,\ldots,j_k$, is proper
and has determinant $1\pm b_n\prod a_{j_ij_{i+1}}\ne1$. This
contradiction finishes the proof.
\end{proof}

The following theorem is not new; the equivalence
$a)\Leftrightarrow b)$ is a particular case
of~\cite[Th.~6]{dobr01}, while $b)\Leftrightarrow c)$ was
discussed (in a more general case of quasitoric manifolds over
arbitrary polytopes) in~\cite{masu99}
and~\cite[Prop.~5.53]{bu-pa02}. We give a proof because it is
needed in the next sections.

\begin{theo}\label{qtcub}
Let $M$ be a quasitoric manifold over a cube, and $\L_\star$ the
corresponding reduced submatrix. The following conditions are
equivalent:
\begin{itemize}
\item[a)] $M$ is equivalent to a Bott tower;

\item[b)] all principal minors of $-\L_\star$ are 1;

\item[c)] $M$ has a $\mathbb T^n$-equivariant almost complex
structure (with the associated omniorientation).
\end{itemize}
\end{theo}
\begin{proof} The implication $\text{b)}\Rightarrow\text{a)}$ follows from Lemma~\ref{minor} and Proposition~\ref{qtbtc}.
The implication $\text{a)}\Rightarrow\text{c)}$ is obvious. Let us
prove $\text{c)}\Rightarrow\text{b)}$. First, we recall
from~\cite[\S4]{masu99}, \cite{dobr01} and~\cite{pano01} the
definition of the \emph{sign} $\sigma(p)$ for a fixed point of the
$\mathbb T^n$-action on $M$. Every fixed point $p$ can be obtained
as the intersection $M_{j_1}\cap\ldots\cap M_{j_n}$ of $n$
characteristic submanifolds, and corresponds to the vertex of $P$
obtained as the intersection $F_{j_1}\cap\ldots\cap F_{j_n}$ of
the corresponding facets. The tangent space to $M$ at $p$
therefore decomposes into the sum of normal subspaces to $M_{j_k}$
for $1\le k\le n$:
\begin{equation}\label{2orie}
  \tau_p(M)={\nu_{j_1}|}_p\oplus\ldots\oplus{\nu_{j_n}|}_p.
\end{equation}
The omniorientation of $M$ provides two different ways of
orienting the above space; we set $\sigma(p)=1$ if these two
orientations coincide and $\sigma(p)=-1$ otherwise. This sign can
be calculated in terms of $P$ and characteristic matrix $\L$ as
\begin{equation}\label{sigma}
  \sigma(p)=\mathop{\mathrm{sign}}\bigl(
  \det(\lambda_{j_1},\ldots,\lambda_{j_n})\cdot
  \det(a_{j_1},\ldots,a_{j_n})\bigr)
\end{equation}
(see~\cite[\S1]{pano01}), where $a_{i}$ denotes a normal vector to
the facet $F_i$ pointing inside the polytope. Now, if $P=\mathbb
I^n$, then every fixed point $p$ corresponds to a vertex given as
\[
  F_{i_1}\cap\ldots\cap F_{i_k}\cap
  F_{n+l_1}\cap\ldots\cap F_{n+l_{n-k}}
\]
for some $1\le i_1<\ldots<i_k\le n$ and $1\le
l_1<\ldots<l_{n-k}\le n$, and we may choose $a_i=e_i$ (the $i$th
basis vector) for $1\le i\le n$ and $a_j=-e_j$ for $n+1\le j\le
2n$. Thus, the expression in the right hand side of~\eqref{sigma}
equals the principal minor of $-\L_\star$ formed by the columns
and rows with numbers $l_1,\ldots,l_{n-k}$. It remains to notice
that in the almost complex case the two orientations
in~\eqref{2orie} coincide, so the sign of every fixed point is~1.
\end{proof}

Recall that a \emph{crosspolytope} is a regular polytope dual to
the cube (in particular, a 3-dimensional crosspolytope is an
octahedron).

\begin{coro}\label{torbt}
Let $X$ be a toric manifold whose associated fan is
combinatorially equivalent to the fan consisting of the cones over
the faces of a crosspolytope. Then $X$ is a Bott tower.
\end{coro}
\begin{proof}
The reduced matrix $\L_\star$ of $X$ is $n\times n$ and all
principal minors of $-\L_\star$ are unit by the same reason as in
the proof of Theorem~\ref{qtcub}. By Lemma~\ref{minor},
$-\L_\star$ can be conjugated to a unipotent uppper triangular
matrix, so $\L$ is of the same form as the characteristic matrix
of a Bott tower. In the toric manifold case, the columns of $\L$
are the primitive vectors along the edges of the fan, so the
combinatorial type of the fan and $\L$ completely determine the
fan. It follows that the fan of $X$ is the same as the fan of a
Bott tower, which implies that $X$ is a Bott tower by the
one-to-one correspondence between fans and toric varieties.
\end{proof}

A toric manifold over a cube satisfies the assumption of
Corollary~\ref{torbt}. It follows that the class of Bott towers
coincides with the class of toric manifolds over a cube, and the
first inclusion in~\eqref{hiera} is an identity. Like
Lemma~\ref{minor}, Corollary~\ref{torbt} is a particular case of a
more general result of Dobrinskaya~\cite[Cor.~7]{dobr01}, which
gives a criterion for a quasitoric manifold over a product of
simplices to be decomposable into a tower of fibre bundles.

\section{Semifree circle actions}\label{semif}
An action of a group on a topological space is called
\emph{semifree} if it is free on the complement to fixed points.
We first show (Theorem~\ref{qtbt} below) that if the torus
$\mathbb T^n$ acting on a quasitoric manifold $M$ over a cube has
a circle subgroup acting semifreely and with isolated fixed
points, then $M$ is a Bott tower. Then we prove that all these
Bott towers are $\mathbb S^1$-equivariantly homeomorphic to a
product of 2-dimensional spheres (with the diagonal $\mathbb
S^1$-action).

A complex $n$-dimensional representation of $\mathbb S^1$ is
determined by a set of weights $k_j\in\Z$ for  $1\le j\le n$. In
appropriate coordinates, an element $s=e^{2\pi i\f}\in \mathbb
S^1$ acts as
\begin{equation}\label{cirac}
  s\cdot(z_1,\ldots,z_n)=(e^{2\pi ik_1\f}z_1,\ldots,e^{2\pi
  ik_n\f}z_n).
\end{equation}
The following statement is straightforward.

\begin{prop}
A representation of $\mathbb S^1$ in $\C^n$ is semifree if and
only if $k_j=\pm1$ for $1\le j\le n$.
\end{prop}

A circle subgroup in $\mathbb T^n$ is determined by an integer
primitive vector $\nu=(\nu_1,\ldots,\nu_n)$:
\begin{equation}\label{cirnu}
  S_\nu^1=\{(e^{2\pi i\nu_1\f},\ldots,e^{2\pi i\nu_n\f})\}
  \subseteq \mathbb T^n.
\end{equation}
We shall consider the tangential representations of $\mathbb T^n$
and its circle subgroups at fixed points of~$M$. The
representation of $\mathbb T^n$ in the tangent space $\tau_p(M)$
at a fixed point $p=M_{j_1}\cap\ldots\cap M_{j_n}$ decomposes into
the sum of non-trivial real two-dimensional representations in the
normal subspaces $\nu_{j_k}$ of $M_{j_k}$. The omniorientation
endows each $\nu_{j_k}$ with a complex structure, therefore
identifying it with $\mathbb C$ and $\tau_pM$ with $\mathbb C^n$.
In these coordinates, the weights of the representation of circle
subgroup~\eqref{cirnu} in $\tau_pM$ can be identified with the
coefficients $k_i=k_i(\nu,p)$, $1\le i\le n$, of the decomposition
of $\nu$ in terms of $\lambda_{j_1},\ldots,\lambda_{j_n}$:
\begin{equation}\label{nucoe}
  \nu=k_1(\nu,p)\lambda_{j_1}+\ldots+k_n(\nu,p)\lambda_{j_n}
\end{equation}
(see, e.g.~\cite[Lemma~2.3]{pano01}).

\begin{coro}\label{sfcoe}
A subcircle $S^1_\nu\subseteq \mathbb T^n$ acts on a quasitoric
manifold $M$ semifreely and with isolated fixed points if and only
if for every vertex $p=F_{j_1}\cap\ldots\cap F_{j_n}$ the
coefficients in~\eqref{nucoe} satisfy $k_i(\nu,p)=\pm1$ for $1\le
i\le n$.
\end{coro}

\begin{theo}\label{qtbt}
Let $M$ be a quasitoric manifold over a cube with reduced
submatrix $\L_\star$. Assume that $M$ admits a semifree circle
subgroup with isolated fixed points. Then $M$ is equivalent to a
Bott tower.
\end{theo}
\begin{proof}
We may assume by induction that every characteristic submanifold
is a Bott tower, so the determinant of every proper principal
minor of $-\L_\star$ is~$1$. Therefore, we are in the situation of
Lemma~\ref{minor}, and $-\L_\star$ is of one of the two types
described there. The second type is ruled out because of the
semifree assumption. Indeed, assume that $\L=\left(E\;|\;B\right)$
where $B$ is matrix~\eqref{exce}, and $S^1_\nu\subseteq \mathbb
T^n$ acts semifreely with isolated fixed points. Applying the
criterion from Corollary~\ref{sfcoe} to the vertex
$p=F_1\cap\ldots\cap F_n$, we obtain $\nu_i=\pm 1$ for $1\le i\le
n$. Now apply the same criterion to the vertex
$p'=F_{n+1}\cap\ldots\cap F_{2n}$. Since the submatrix formed by
the corresponding columns of $\L$ is exactly $B$, we have $\det
B=\pm1$. This implies that at least one of $b_i$ equals $\pm1$,
i.e. at least one of the rows of $B$ consists only of two $\pm1$'s
and zeros. Therefore, if all the coefficients $k_i(\nu,p')$ in the
expression
$\nu=k_1(\nu,p')\lambda_{n+1}+\ldots+k_n(\nu,p')\lambda_{2n}$ are
$\pm1$, then at least one component $\nu_i$ in the left hand side
is even. A contradiction.
\end{proof}

Our next result shows that a Bott tower with a semifree circle
subgroup and isolated fixed points is topologically (or even
$\mathbb S^1$-equivariantly) trivial, i.e. homeomorphic to a
product of 2-spheres. Let $t$ (respectively, $\C$) denote the
standard (respectively, trivial) complex one-dimensional $\mathbb
S^1$-\-re\-pre\-sen\-ta\-tion. A product bundle with fibre $V$
will be denoted by $\underline{V}$. We say that an action of a
group $G$ on a Bott tower $B^{2n}$ \emph{preserves the tower
structure} if for each stage
$B^{2k}=P(\underline{\C}\oplus\xi_{k-1})$, \ $k\le n$, the line
bundle $\xi_{k-1}$ is $G$-equivariant. The intrinsic $\mathbb
T^n$-action obviously preserves the tower structure.

\begin{theo}\label{btpro}
Assume that a Bott tower $B^{2n}$ admits a semifree $\mathbb
S^1$-action with isolated fixed points preserving the tower
structure. Then $B^{2n}$ is $\mathbb S^1$-equivariantly
diffeomorphic to the product $(P(\C\oplus t))^n$.
\end{theo}
\begin{proof}We may
assume by induction that the $(n-1)$-stage of the Bott tower is
$(P(\C\oplus t))^{n-1}$ and $B^{2n}=P(\underline{\C}\oplus \xi)$
for some $\mathbb S^1$-line bundle $\xi$ over ${(P(\C\oplus
t))^{n-1}}$.

Let $\gamma$ be the unique $\mathbb S^1$-line bundle over
$P(\C\oplus t)\cong\C P^1$ whose underlying bundle is the
canonical line bundle and
\begin{equation} \label{gamma}
\gamma|_{(1:0)}=\C \qquad\text{and}\qquad \gamma|_{(0:1)}=t.
\end{equation}
Denote by $x\in{H^2(P(\C\oplus t))}$ the first Chern class of
$\gamma$, and let $x_i\in H^2(P(\C\oplus t)^{n-1})$ be the
pullback of $x$ by the projection $\pi_i$ onto the $i$-th factor.
The first Chern class of $\xi$ may then be written as
$\sum_{i=1}^{n-1}a_i x_i$ with $a_i\in \Z$. The $\mathbb S^1$-line
bundles $\xi$ and $\otimes_{i=1}^{n-1}\pi_i^*(\gamma^{a_i})$ have
the same underlying bundles; so there is an integer $k$ such that
\begin{equation} \label{xi}
  \xi= t^k\bigotimes_{i=1}^{n-1}\pi_i^*(\gamma^{a_i})
\end{equation}
as $\mathbb S^1$-line bundles, see \cite[Cor. 4.2]{ha-yo76}.

We encode the fixed points of $P(\C\oplus t)^{n-1}$ as
$(p_1^{\epsilon_1},\dots,p_{n-1}^{\epsilon_{n-1}})$ where
$\epsilon_i=0$ or $1$ and $p_i^{\epsilon_i}$ denotes $(1:0)$ if
$\epsilon_i=0$ and $(0:1)$ if $\epsilon_i=1$.  It follows
from~\eqref{gamma} and~\eqref{xi} that
\[
  \xi|_{(p_1^{\epsilon_1},\dots, p_{n-1}^{\epsilon_{n-1}})}
  =t^{k+\sum_{i=1}^{n-1} \epsilon_i a_i}.
\]
The $\mathbb S^1$-action on $B^{2n}=P(\underline{\C}\oplus\xi)$ is
semifree if and only if $|k+\sum_{i=1}^{n-1} \epsilon_i a_i|=1$
for all possible values of $\epsilon_i$'s. Setting $\epsilon_i=0$
for all $i$ we get $|k|=1$. Assume $k=1$ (the case $k=-1$ is
treated similarly). Then $(a_1,\dots,a_{n-1})=(0,\dots,0)$ or
$(0,\dots,0,-2,0,\dots,0)$. In the former case,
$\xi=\underline{t}$ and $B^{2n}=P(\underline{\C}\oplus\xi)\cong
P(\C\oplus t)^n$. In the latter case we have
$\xi=t\pi_i^*(\gamma^{-2})$ for some $i$, so that
$B^{2n}=P(\underline{\C}\oplus\xi)$ is the pullback of
$P(\underline{\C}\oplus t\gamma^{-2})$ by the projection $\pi_i$.
Since for any $\mathbb S^1$-vector bundle $E$ and $\mathbb
S^1$-line bundle $\eta$ the projectivisations $P(E)$ and
$P(E\otimes\eta)$ are $\mathbb S^1$-diffeomorphic, we have
$P(\underline{\C}\oplus t\gamma^{-2})\cong P(\gamma\oplus
t\gamma^{-1})$. The first Chern class of $\gamma\oplus
t\gamma^{-1}$ is zero, so its underlying bundle is trivial. The
$\mathbb S^1$-representation in the fibre of $\gamma\oplus
t\gamma^{-1}$ over a fixed point is $\C\oplus t$ by~\eqref{gamma}.
Therefore, $\gamma\oplus t\gamma^{-1}=\underline{\mathbb
C}\oplus\underline t$ as $\mathbb S^1$-bundles. It follows that
$P(\underline{\C}\oplus t\gamma^{-2})\cong P(\underline{\mathbb
C}\oplus\underline t)$, which finishes the proof.
\end{proof}

\begin{rema}
The diffeomorphism of Theorem~\ref{btpro} is not a diffeomorphism
of $\mathbb T^n$-manifolds.
\end{rema}

Our next aim is to generalise Theorem~\ref{btpro} to quasitoric
manifolds. Although the result does not hold for all quasitoric
manifolds (see the next Example), it remains true if we
additionally require the quotient polytope to be a cube.

\begin{exam}
Let $M$ be the 4-dimensional quasitoric manifold over a $2k$-gon
with characteristic matrix of the form
$\begin{pmatrix}1&0&1&0&\ldots&1&0\\0&1&0&1&\ldots&0&1\end{pmatrix}$.
Corollary~\ref{sfcoe} shows that the circle subgroup determined by
the vector $\nu=(1,1)$ acts semifreely on $M$, but the quotient of
$M$ is not a 2-cube if $k>2$, so $M$ can not be homeomorphic to
the product of spheres (it may be shown that $M$ is the connected
sum of $k-1$ copies of $S^2\times S^2$).
\end{exam}

Surprisingly, quasitoric manifolds over polygons provide the only
essential source of counterexamples, see Theorem~\ref{qtmsf} below
for the precise statement.

\begin{lemm}\label{spcub}
A simple polytope $P$ of dimension $n\ge2$ all of whose 2-faces
are 4-gons is combinatorially equivalent to a cube.
\end{lemm}
\begin{proof}
We may assume by induction that all facets of $P$ are cubes, and
prove that $P$ is a cube. The statement can be found
in~\cite[Exercise~0.1]{zieg95}, but we include the proof for the
sake of completeness. We shall prove the dual statement about
simplicial polytopes. The simplicial polytope dual to $P$ is a
crosspolytope; we call its boundary $K$ (which is a sphere
triangulation) a \emph{cross-complex}. Recall that the \emph{star}
of a vertex $v$ in a simplicial complex $K$ is the subcomplex $\st
v$ consisting of all simplices that together with $v$ span a
simplex, and \emph{link} of $v$ is the subcomplex $\lk v\subset\st
v$ consisting of simplices not containing~$v$. The duality between
$P$ and $K$ extends to the duality between the facets of $P$
(which are $(n-1)$-dimensional simple polytopes) and the links of
vertices of $K$ (which are triangulations of $(n-2)$-spheres). The
dual statement follows from the next lemma (which is slightly more
general, as we do not assume $K$ to be the boundary of a
simplicial polytope there).
\end{proof}

\begin{lemm}\label{cross}
Let $K$ be a connected simplicial complex of dimension $k\ge 2$.
If the link of each vertex of $K$ is a crosscomplex of dimension
$k-1$, then $K$ is a crosscomplex.
\end{lemm}
\begin{proof}
Let $v$ be a vertex of $K$. By the assumption, $\lk v$ is a
crosscomplex of dimension $k-1$, so every vertex $p\in\lk v$ has a
unique vertex $q\in\lk v$ which is not joined to $p$ by an edge in
$\lk v$. Still, $p$ and $q$ may be joined by an edge in~$K$, so we
consider the two cases.

{\sc Case 1.} Suppose that there is a pair of vertices $p,q$ in
$\lk v$ that are not joined by an edge in $K$. Let $\mathcal R$ be
the set of other vertices of $\lk v$. The cardinality of $\mathcal
R$ is $2(k-1)$. The link $\lk p$ is a crosspolytope, so it has
$2k$ vertices, and contains $v$ and all elements in $\mathcal R$.
Since $q$ is not joined to $p$ by an edge in $K$, $q$ is not in
$\lk p$; so there is another vertex $p'\in\lk p$, \ $p'\notin
v\cup\mathcal R$. Similarly, we have $q'\in\lk q$, \ $q'\notin
v\cup\mathcal R$. Now take any $r\in\mathcal R$ and consider $\lk
r$. Since $\lk v$ is a $(k-1)$-dimensional crosscomplex, $r$ is
joined to $2(k-1)$ vertices of $\lk v$ by edges in $\lk v$. We
also know that $r$ is joined to $v$, $p'$ and $q'$. But since $\lk
r$ is also a $(k-1)$-dimensional crosscomplex, $r$ may be joined
to only $2k$ vertices. Therefore $p'=q'$, which implies that $K$
is a crosscomplex.

{\sc Case 2.} Suppose that every pair of vertices in $\lk v$ is
joined by an edge in $K$. This will lead us to a contradiction.
Each vertex $u$ in $\lk v$ is joined to $v$ and all vertices in
$\lk v$ except $u$ itself. These are all vertices joined to $u$ by
edges, because the number of vertices in $\lk u$ is $2k$. This
means that any pair of vertices in $K$ is joined by an edge, and
that $K$ has exactly $2k+1$ vertices. The number of $k$-simplices
meeting at each vertex is $2^k$, and a $k$-simplex has $k+1$
vertices. So the total number of $k$-simplices in $K$ is
$2^k(2k+1)/(k+1)$.

Now we calculate the total number of $k$-simplices in $K$ in a
different way. Let $\sigma$ be a $k$-simplex of $K$ not containing
$v$. Then $\sigma$ contains a pair of vertices, say $p$ and $q$,
that are not joined by an edge in $\lk v$ (otherwise $\sigma$
itself must be in $\lk v$, since $\lk v$ is a crosscomplex). Let
$L$ be the link of $p$ in $\lk v$. Then $L$ is a crosscomplex of
dimension $k-2$, and it also coincides with the link of $q$ in
$\lk v$. We have that $\lk p$ is the join $L*\{ v,q\}$, because
both subcomplexes have the same vertex sets and both are
crosscomplexes, and similarly $\lk q=L*\{ v,p\}$. Since $\sigma$
contains $p$ and $q$, it follows that $\sigma=\tau*p*q$ for some
$(k-2)$-simplex $\tau\in L$. Therefore, $\sigma$ has at least two
faces of dimension $(k-1)$ in $\lk v$, namely $\tau*p$ and
$\tau*q$. Neither of these can be a face of another $k$-simplex
not containing $v$, because every $(k-1)$-simplex of $K$ is a face
of exactly two $k$-simplices by the assumption, and $\tau*p$ is
also a face of $\tau*p*v$ and $\tau*q$ is also a face of
$\tau*q*v$. It follows that the number of $k$-simplices not
containing $v$ is no more than half of the number of
$(k-1)$-simplices in $\lk v$. The number of $k$ simplices
containing $v$ equals the number of $(k-1)$-simplices in $\lk v$.
The latter number is $2^k$, so the total number of $k$-simplices
in $K$ is $\le 2^{k-1}+2^k$, which is less than $2^k(2k+1)/(k+1)$
if $k\ge 2$. This contradiction finishes the proof.
\end{proof}

\begin{rema}
Another proof of Lemma~\ref{spcub} may be given by establishing a
non-degenerate simplicial map from $K$ to a crosscomplex. Such a
map is a topological (non-ramified) covering of a sphere by a
sphere, so it must be an isomorphism for $n\ge3$. This approach
was used in~\cite{ilin06}.
\end{rema}

\begin{theo}\label{qtmsf}
Assume that a quasitoric manifold $M$ admits a semifreely acting
subcircle with isolated fixed points, and every 2-face of the
quotient polytope $P$ is a 4-gon. Then $M$ is $\mathbb
S^1$-equivariantly homeomorphic to a product of 2-dimensional
spheres.
\end{theo}
\begin{proof}
By Lemma~\ref{spcub}, the orbit polytope is a cube. By
Theorem~\ref{qtbt}, $M$ is equivalent to a Bott tower. Finally, by
Theorem~\ref{btpro} it is $\mathbb S^1$-homeomorphic to a product
of spheres.
\end{proof}

We are also able to deduce the main result of
Ilinskii~\cite{ilin06}:

\begin{coro}\label{ilins}
A toric manifold $X$ with a circle subgroup acting semifreely and
with isolated fixed points is diffeomorphic to a product of
2-spheres.
\end{coro}
\begin{proof}
By Theorem~\ref{btpro}, it suffices to show that $X$ is a Bott
tower. To do this we apply Corollary~\ref{torbt}, that is, we show
that the fan corresponding to $X$ is combinatorially equivalent to
the fan over a crosspolytope. The semifree circle subgroup acting
on $X$ also acts semifreely and with isolated fixed points on
every characteristic submanifold $X_j$ of $X$. Using induction by
dimension and Lemma~\ref{cross} we reduce the statement to the
complex 2-dimensional case, that is, we need to show that the
quotient polytope of a complex 2-dimensional toric manifold with a
semifree circle subgroup action and isolated fixed points is a
4-gon. (Note that a complex 2-dimensional toric manifold is always
projective, so that we can work with polytopes instead of fans.)
The following case-by-case analysis is just a reformulation of the
argument from~\cite[\S3]{ilin06}; we give it here for the sake of
completeness.

Let $\Sigma$ be the fan corresponding to our complex 2-dimensional
toric manifold. The one dimensional cones of $\Sigma$ correspond
to the facets (or edges) of the quotient polygon $P^2$. We need to
show that there are exactly 4 of them. The values of the
characteristic function on the facets of $P^2$ are given by the
primitive vectors generating the corresponding one-dimensional
cones of $\Sigma$. Let $\nu$ be the vector generating the semifree
circle subgroup. We may choose the initial vertex $p$ of $P^2$ so
that $\nu$ belongs to the 2-dimensional cone of $\Sigma$
corresponding to~$p$. Then we index the primitive vectors
$\lambda_i$, $1\le i\le m$, so that $\nu$ is in the cone generated
by $\lambda_1$ and $\lambda_2$, and any two consecutive vectors
span a two-dimensional cone (see Fig.~\ref{figur}).
\begin{figure}[h]
\begin{center}
\begin{picture}(60,40)
  \put(40,20){\vector(1,0){20}}
  \put(40,20){\vector(1,1){20}}
  \put(40,20){\vector(0,1){20}}
  \put(40,20){\vector(-1,0){20}}
  \put(40,20){\vector(-2,-1){40}}
  \put(59.5,16){$\lambda_1$}
  \put(59.5,37){$\nu$}
  \put(36,37){$\lambda_2$}
  \put(19,21.5){$\lambda_3$}
  \put(0,3.5){$\lambda_m$}
  \multiput(9,9)(1.5,1.5){3}{$\cdot$}
\end{picture}
\caption{ } \label{figur}
\end{center}
\end{figure}
This provides us with a refined $2\times m$ characteristic matrix
$\L$. We have $\lambda_1=(1,0)$ and $\lambda_2=(0,1)$, and
applying the criterion from Corollary~\ref{sfcoe} to the first
cone $\langle\lambda_1,\lambda_2\rangle$ (i.e., to the initial
vertex of the polygon), we obtain $\nu=(1,1)$.

Now consider the second cone. The non-singularity condition gives
us $\det(\lambda_2,\lambda_3)=1$, whence $\lambda_3=(-1,*)$.
Writing $\nu=k_1\lambda_2+k_2\lambda_3$ and applying
Corollary~\ref{sfcoe} to the second cone
$\langle\lambda_2,\lambda_3\rangle$, we obtain
\[
  (1,1)=\pm(0,1)\pm(-1,*).
\]
Therefore, $\lambda_3=(-1,0)$ or $\lambda_3=(-1,-2)$. Similarly,
considering the last cone $\langle\lambda_m,\lambda_1\rangle$, we
get $\lambda_m=(*,-1)$, and then, applying Corollary~\ref{sfcoe},
$\lambda_m=(0,-1)$ or $\lambda_m=(-2,-1)$. The case
$\lambda_3=(-1,-2)$ and $\lambda_m=(-2,-1)$ is impossible as then
the second and the last cones overlap.

Assume $\lambda_3=(-1,0)$. Then a similar analysis applied to the
third cone $\langle\lambda_3,\lambda_4\rangle$ shows that
$\lambda_4=(0,-1)$ or $\lambda_4=(-2,-1)$. Therefore,
$\lambda_4=\lambda_m$ (otherwise cones overlap).

Similarly, if $\lambda_m=(0,-1)$, we get $\lambda_{m-1}=(-1,0)$ or
$\lambda_{m-1}=(-1,-2)$. Therefore, $\lambda_{m-1}=\lambda_3$.

In any case, we have $m=4$ and $P^2$ is a 4-gon. This completes
the proof.
\end{proof}

Note that the proof above leaves three possibilities for the
vectors $\lambda_3$ and $\lambda_4$ of the corresponding
2-dimensional fan: $(-1,0)$ and $(0,-1)$, or $(-1,0)$ and
$(-2,-1)$, or $(-1,-2)$ and $(0,-1)$, and the last two pairs are
equivalent by an orientation reversing automorphism of $T^2$. The
corresponding reduced submatrices are
$\begin{pmatrix}-1&0\\0&-1\end{pmatrix}$ and
$\begin{pmatrix}-1&-2\\0&-1\end{pmatrix}$. The first one
corresponds to $\C P^1\times\C P^1$, and the second to a
non-trivial Bott tower (Hirzebruch surface) with $a_{12}=-2$. We
finish this section by determining the class of
matrices~\eqref{amatr} corresponding to our specific class of Bott
towers explicitly, for arbitrary dimension.

\begin{theo}\label{sfdec}
A Bott tower $B^{2n}$ admits a semifree circle subgroup with
isolated fixed points if and only if its matrix~\eqref{amatr}
satisfies the identity
\[
  \frac12(E-A)=C_1C_2\cdots C_n
\]
where $C_k$ is either a unit matrix or a unipotent upper
triangular matrix with only one non-zero element $c_{i_kk}=1$ in
the $k$-th column off the diagonal, for $1\le k\le n$.
\end{theo}
\begin{proof}
First assume that $B^{2n}$ admits a semifree circle subgroup with
isolated fixed points. We have two sets of multiplicative
generators for $H^*(B^{2n})$: the set $\{u_1,\ldots,u_n\}$ from
Lemma~\ref{cohbt}, satisfying~\eqref{crobt}, and the set
$\{x_1,\ldots,x_n\}$, satisfying $x_i^2=0$, $1\le i\le n$. The
reduced sets with $i\le k$ can be considered as the corresponding
sets of generators for the $k$-th stage $B^{2k}$. As it is clear
from the proof of Theorem~\ref{btpro}, we have
$c_1(\xi_{k-1})=-2c_{i_kk}x_{i_k}$ for some $i_k<k$, where
$c_{i_kk}=1$ or $0$. From $u_k^2+2c_{i_kk}x_{i_k}u_k=0$ we get
$x_k=u_k+c_{i_kk}x_{i_k}$. In other words, the transition matrix
$C_k$ from the basis $x_1,\ldots,x_{k-1},u_k,\ldots,u_n$ of
$H^2(B^{2n})$ to $x_1,\ldots,x_k,u_{k+1},\ldots,u_n$ may have only
one non-zero entry off the diagonal, which is $c_{i_kk}$. The
transition matrix from $u_1,\ldots,u_n$ to $x_1,\ldots,x_n$ is
therefore the product $D=C_1C_2\cdots C_n$ (note that $C_1$ is the
unit matrix as $x_1=u_1$). Then $D=(d_{jk})$ is a unipotent upper
triangular matrix consisting of zeros and units, we have
$x_k=\sum_{j=1}^nd_{jk}u_j$, and
\[\textstyle
  0=x_k^2=\bigl(u_k+\sum_{j=1}^{k-1}d_{jk}u_j\bigr)^2=u_k^2+
  2\sum_{j=1}^{k-1}d_{jk}u_ju_k+\ldots, \quad
  k=1,\ldots,n.
\]
On the other hand, $0=u_k^2-\sum_{j=1}^{k-1}a_{jk}u_ju_k$
by~\eqref{crobt}. Comparing the coefficients of $u_ju_k$ for $1\le
j\le k-1$ in the last two equations and noting that these elements
form a basis of the quotient space $H^4(B^{2k})/H^4(B^{2(k-1)})$,
we obtain $2d_{jk}=-a_{jk}$ for $1\le j<k\le n$. As both $D$ and
$-A$ are unipotent upper triangular matrices, this implies
$2D=E-A$, as needed.

On the other hand, if $A$ satisfies $E-A=2C_1C_2\cdots C_n$, then
for the corresponding Bott tower we have
$\xi_{k-1}=\pi_{i_k}^*(\gamma^{-2c_{i_kk}})$. Therefore we may
choose a circle subgroup in such a way that $\xi_{k-1}$ becomes
$t\pi_{i_k}^*(\gamma^{-2c_{i_kk}})$ as an $\mathbb S^1$-bundle,
for $1<k\le n$. This circle subgroup acts semifreely and with
isolated fixed points by the same argument as in the proof of
Theorem~\ref{btpro}.
\end{proof}

\begin{exam}
It follows from Theorem~\ref{sfdec} that if a Bott tower admits a
semifree circle subgroup with isolated fixed points, then
matrix~\eqref{amatr} may have only 0 or $-2$ entries above the
diagonal. However, the condition of Theorem~\ref{sfdec} is
stronger. For $n=3$, the matrix
\[
  A=\begin{pmatrix} -1&0&-2\\ 0&-1&-2\\ 0&0&-1 \end{pmatrix}
\]
does not occur for our class of Bott towers (since $(E-A)/2$
cannot be decomposed as $C_1C_2C_3$), while any other matrix with
0 or $-2$ above the diagonal occurs. For example, if
\[
  A=\begin{pmatrix} -1&-2&-2\\ 0&-1&0\\ 0&0&-1 \end{pmatrix},
\]
then we have
\[
  \frac12(E-A)=\begin{pmatrix} 1&1&1\\ 0&1&0\\ 0&0&1 \end{pmatrix}=
  \begin{pmatrix} 1&0&0\\ 0&1&0\\ 0&0&1 \end{pmatrix}
  \begin{pmatrix} 1&1&0\\ 0&1&0\\ 0&0&1 \end{pmatrix}
  \begin{pmatrix} 1&0&1\\ 0&1&0\\ 0&0&1 \end{pmatrix}.
\]
\end{exam}

It is clear that not every Bott tower homeomorphic to a product of
spheres admits a semifree subcircle action with isolated fixed
points (the latter condition is stronger even for $n=2$). We shall
consider the class of Bott towers that are homeomorphic to a
product of spheres in the next section.

\section{Topological classification and cohomology}\label{tocla}
The following statement shows that Bott towers diffeomorphic to
products of spheres can be detected by their cohomology rings,
thereby providing a partial answer to Problem~\ref{cohcl}:

\begin{theo}\label{trivc}
A Bott tower $B^{2n}$ is diffeomorphic to $(\C P^1)^n$ if
$H^*(B^{2n})\cong H^*((\C P^1)^n)$ as rings.
\end{theo}
\begin{proof}
From Lemma~\ref{cohbt} we get
\[
  H^*(B^{2n})=H^*\bigl(B^{2n-2}\bigr)[u_n]\bigl/
  \bigl(u_n^2-c_1(\xi_{n-1})u_n\bigr).
\]
We may therefore write any element of $H^2(B^{2n})$ as $x+bu_n$
where $x\in H^2(B^{2n-2})$ and $b\in\Z$. Since
\[
  (x+bu_n)^2=x^2+2bxu_n+b^2u_n^2=x^2+b(2x+bc_1(\xi_{n-1}))u_n,
\]
the square of $x+bu_n$ with $b\not=0$ is zero if and only if
$x^2=0$ and $2x+bc_1(\xi_{n-1})=0$. This shows that the elements
$x+bu_n$ with $b\ne0$ whose squares are zero span a rank-one
subgroup in $H^2(B^{2n})$.

By assumption, there is a basis $\{x_1,\dots,x_n\}$ of
$H^2(B^{2n})$ such that $x_i^2=0$ for all $i$. By the observation
above, we may assume that $x_1,\dots,x_{n-1}$ are in
$H^2(B^{2n-2})$. Because $\{x_1,\dots,x_n\}$ is a basis, $x_n$ is
not in $H^2(B^{2n-2})$ and we may assume
$x_n=\sum_{i=1}^{n-1}b_ix_i+u_n$ with some $b_i\in\Z$. A product
of the form $\prod_{i\in I}x_i$, where $I$ is a subset of
$\{1,\dots,n\}$, belongs to $H^*(B^{2n-2})$ if and only if
$n\notin I$. This shows that $H^*(B^{2n-2})$ is generated by
$x_1,\dots,x_{n-1}$ and has the same cohomology ring as $(\C
P^1)^{n-1}$. Therefore, we may assume $B^{2n-2}\cong (\C
P^1)^{n-1}$ by induction.

Writing $c_1(\xi_{n-1})=\sum_{i=1}^{n-1}a_ix_i$, we see that
\[\textstyle
  0=x_n^2=\bigl(u_n+\sum_{i=1}^{n-1} b_ix_i\bigr)^2
  =\sum_{i=1}^{n-1}(a_i+2b_i)x_iu_n+
  \bigl(\sum_{i=1}^{n-1} b_ix_i\bigr)^2.
\]
This may hold only if at most one $a_i$ is non-zero because
$x_ix_j$ $(i<j)$ and $x_iu_n$ form an additive basis of
$H^4(B^{2n})$. Therefore, $\xi_{n-1}$ is the pullback of
$\gamma^{-2b_i}$ over $\C P^1$ by a projection $B^{2n-2}=(\C
P^1)^{n-1}\to \C P^1$. Since
$P(\underline{\C}\oplus\gamma^{-2b_i})$ is a product bundle (see
Example~\ref{2dbto}), so is
$B^{2n}=P(\underline{\C}\oplus\xi_{n-1})$.
\end{proof}

We can now also effectively describe the class of
matrices~\eqref{amatr} corresponding to Bott towers which are
diffeomorphic to a product of 2-spheres.

\begin{theo}
A Bott tower $B^{2n}$ is diffeomorphic to $(\C P^1)^n$ if and only
if the corresponding matrix~\eqref{amatr} satisfies the identity
\[
  \frac12(E-A)=C_1C_2\cdots C_n
\]
where each $C_k$ is a unipotent upper triangular matrix which may
have only one non-zero element in the $k$-th column off the
diagonal, for $1\le k\le n$.
\end{theo}
\begin{proof}
We use the same argument as in the proof of Theorem~\ref{sfdec}.
The only difference is that $c_{i_kk}$ in
$c_1(\xi_{k-1})=-2c_{i_kk}x_{i_k}$ may now be an arbitrary
integer.
\end{proof}

In the rest of this section we generalise the result of
Theorem~\ref{trivc} to arbitrary quasitoric manifold, but only in
the topological category (see Theorem~\ref{qtcpn}).

We start by analysing the algebraic structure of the cohomology of
quasitoric manifolds over a cube. Although it is possible to make
this analysis over $\Z$, it is more convenient for our purposes to
reduce the coefficients modulo~2. Let $S$ be a graded algebra over
$\Z/2$ generated by degree one elements $x_1,\ldots,x_n$. We refer
to $S$ as a \emph{Bott quadratic algebra} (simply
\emph{BQ-algebra}) of rank $n$ if it satisfies the following two
two properties:
\smallskip
\begin{enumerate}
\item[(P1)] $x_k^2=\sum_{i<k}a_{ik} x_ix_k$ \ with $a_{ik}\in\Z/2$
for $1\le k\le n$. (In particular, $x_1^2=0$.)

\item[(P2)] $\prod_{i=1}^n x_i\not=0$.
\end{enumerate}

If $B^{2n}$ is a Bott tower, then~\eqref{crobt} implies that
$H^*(B^{2n};\Z/2)$ is a $\BQ$-algebra with the degree doubled,
which explains our terminology. Our arguments below work for a
wider class of algebras with~(P1) replaced by the weaker property:
\begin{enumerate}
\item[(P1')] $x_k^2=\sum_{i<j\le k}a_{ijk} x_ix_j$ \ with
$a_{ijk}\in\Z/2$ for $1\le k\le n$.
\end{enumerate}

Because of (P1) we can express any element of $S$ as a linear
combination of square-free monomials. We denote such a monomial
$x_{i_1}\ldots x_{i_s}$ by $x_I$, where $I=\{i_1,\dots,i_s\}$.

\begin{lemm}\label{basis}
The elements $x_I$ for all subsets $I\in\{1,\dots,n\}$ form an
additive basis of~$S$. In particular, $\dim S_q=\binom{n}{q}$
where $S_q$ denotes the graded component of degree~$q$.
\end{lemm}
\begin{proof} (P1) implies that the set $\{x_I\}$
additively generates $S$. We order monomials $x_I$ using the
inverse lexicographical order on subsets of $\{1,\ldots,n\}$.
Namely, if $I=\{i_1,\dots,i_s\}$ with $i_1<\ldots<i_s$ and
similarly $J=\{j_1,\dots,j_s\}$, then we set $x_I<x_J$ if
$i_k<j_k$ and $i_q=j_q$ for $k+1\le q\le s$.

Suppose that there is a non-trivial linear relation for $x_I$, and
let $x_J$ be the maximal monomial appearing in this relation. Then
we may use the relation to replace the subfactor $x_J$ in
$\prod_{i=1}^nx_i$, and then use~(P1) repeatedly whenever $x_k^2$
appears, until we end up at zero. This contradicts~(P2).
Therefore, there are no non-trivial linear relations among
$x_I$'s.
\end{proof}

\begin{lemm} \label{BQ}
Suppose we have a surjective graded homomorphism $f$ from $S$ to a
graded algebra $S'$ over $\Z/2$ satisfying $S'_{n-1}\not=0$. Then
the dimension of the kernel of $f\colon S_1\to S'_1$ is at most
one. Moreover, if the dimension of the kernel is exactly one, then
$S'$ is a $\BQ$-algebra of rank $n-1$.
\end{lemm}
\begin{proof}
We denote $f(x_i)$ by $\bar x_i$. Then~(P1) holds for $\bar
x_1,\dots,\bar x_n$. Suppose the dimension of the kernel is at
least two. Then there exist $p>q\ge 1$ such that
\begin{equation} \label{ypyq}
\bar x_p=\sum_{i<p}b_i\bar x_i,\qquad \bar x_q=\sum_{j<q}c_j\bar
x_j
\end{equation}
where $b_i,c_j\in\Z/2$. By Lemma~\ref{basis}, $S_{n-1}$ is
generated by the elements $x_I$ with $|I|=n-1$. We shall show that
$\bar x_I=0$ for any such~$I$, which contradicts the assumption
$S'_{n-1}\not=0$.

First assume $q\ge 2$. Since $|I|=n-1$, $I$ contains $p$ or $q$.
We replace $\bar x_p$ and $\bar x_q$ in $\bar x_I$ by~\eqref{ypyq}
and use~(P1) repeatedly whenever $\bar x_k^2$ appears. We end up
at zero.

If $q=1$, i.e., $\bar x_1=0$, then it suffices to see that $\bar
x_I=0$ for $I=\{2,3,\dots,n\}$.  We replace $\bar x_p$ in $\bar
x_I$ by~\eqref{ypyq} and use (P1) repeatedly whenever $\bar x_k^2$
appears for $k\ge 2$. Then each term in the final expression
contains a factor $\bar x_1$, which is zero.

Now we prove the second statement of the lemma. By the assumption,
the elements $\bar x_i$ satisfy a non-trivial linear relation. Let
$\bar x_j$ be the maximal element appearing in this relation. We
can eliminate $\bar x_j$ in $S'$ using this linear relation
and~(P1). Then~(P1) holds for $\bar x_i$'s with $i\not=j$.
Therefore, $S'_{n-1}$ is generated by $\prod_{i\not=j}\bar x_i$.
This element is non-zero because $S'_{n-1}\not=0$. This proves
that $S'$ is a $\BQ$-algebra of rank $n-1$.
\end{proof}

\begin{theo} \label{cube}
Let $M$ be a quasitoric manifold with the quotient polytope~$P$.
Then $H^{2*}(M;\Z/2)$ is a $\BQ$-algebra of rank $n$ if and only
if $P$ is an $n$-cube.
\end{theo}
\begin{proof}
Assume that $P$ is an $n$-cube. Since every principal minor of
$\L_\star$ is 1 modulo 2, we deduce that $\L_\star$ is conjugate
to a unipotent upper triangular matrix by the same argument as in
the proof of Lemma~\ref{minor}. (Matrix~\eqref{exce} in which all
$b_i$ are non-zero modulo 2 for $1\le i\le n$ cannot occur because
its determinant is zero modulo 2.) Then it follows from
from~\eqref{redvi} that $H^{2*}(M;\Z/2)$ is a $\BQ$-algebra of
rank $n$.

Now assume that $H^{2*}(M;\Z/2)$ is a $\BQ$-algebra. Let
$b_{r}(M)$ denote the $r$-th Betti number of $M$, and let $f_s(P)$
denote the number of faces of $P$ of codimension $s+1$. Then
\[
  b_2(M)=f_0(P)-n, \qquad b_4(M)=f_1(P)-(n-1)f_0(P)+\binom{n}{n-2}
\]
(see \cite[Th.~3.1]{da-ja91}), and we obtain from
Lemma~\ref{basis} that
\begin{equation} \label{f0f1}
  f_0(P)=2n, \qquad f_1(P)=2n(n-1).
\end{equation}

{\sloppy For every characteristic submanifold $M_i$ the
restriction map $H^*(M;\Z/2)\to H^*(M_i;\Z/2)$ is
surjective~\cite[Lemma~2.3]{ma-pa06}. It follows from
Lemmas~\ref{basis} and~\ref{BQ} that $b_2(M_i)\ge b_2(M)-1=n-1$.
Therefore,
\begin{equation}\label{f0Qi}
  f_0(F_i)=(n-1)+b_2(M_i)\ge 2(n-1),
\end{equation}
where $F_i$ denotes the facet corresponding to $M_i$, and
\[
  f_1(P)=\frac{1}{2}\sum_{i=1}^{2n}f_0(F_i)\ge 2n(n-1).
\]
Comparing this with~\eqref{f0f1} we obtain that the equality holds
in~\eqref{f0Qi} for every $i$, that is, $b_2(M_i)=n-1$. This
implies that the kernel of $H^2(M;\Z/2)\to H^2(M_i;\Z/2)$ is
one-dimensional, so $H^{2*}(M_i;\Z/2)$ is a $\BQ$-algebra of rank
$n-1$ by Lemma~\ref{BQ}. This enables us to use induction on~$n$.

}
When $n=2$, equations~\eqref{f0f1} imply that $P$ is
combinatorially a square. Suppose that the theorem holds for $n-1$
with $n\ge 3$. Since $H^{2*}(M_i)$ is a $\BQ$-algebra of rank
$n-1$, every facet of $P$ is an $(n-1)$-cube; in particular, every
2-face of $P$ is a square. Then $P$ is an $n$-cube by virtue of
Lemma~\ref{spcub}.
\end{proof}

\begin{lemm}\label{Qcoef}
Let $M$ be a quasitoric manifold over an $n$-cube. If
$H^*(M;\Q)\cong H^*((\C P^1)^n;\Q)$ as rings, then $M$ is
equivalent to a Bott tower.
\end{lemm}
\begin{proof}
By the assumption, there are elements $y_1,\ldots,y_n$ in
$H^2(M;\Q)$ generating $H^*(M;\Q)$ and satisfying $y^2_i=0$ for
$1\le i\le n$. Let $M_i\subset M$ be a characteristic submanifold,
and denote the restriction of $y_k$ to $H^2(M_i;\Q)$ by $\bar
y_k$. Then $\bar y_1,\dots,\bar y_n$ generate $H^*(M_i;\Q)$ as a
ring because $H^*(M;\Q)\to H^*(M_i;\Q)$ is surjective. Since
$b_2(M_i)=n-1$, there is a non-trivial linear relation for the
elements $\bar y_k$. Using this linear relation, we can eliminate
one generator, say $\bar y_n$, and get a surjective map $\Q[\bar
y_1,\dots,\bar y_{n-1}]/(\bar y_1^2,\dots,\bar y_{n-1}^2)\to
H^*(M_i;\Q)$. Since the dimensions of components of degree $2q$ in
both rings equal $\binom{n-1}{q}$, the surjective map is an
isomorphism. Therefore, $H^*(M_i;\Q)\cong H^*((\C P^1)^{n-1};\Q)$,
so we may use an induction argument and assume that every $M_i$ is
a Bott tower.

Let $\L_\star$ be the reduced submatrix of $M$. It follows from
Lemma~\ref{minor} that $-\L_\star$ is conjugate to a unipotent
upper triangular matrix or to matrix~\eqref{exce} with non-zero
entries $b_i$ for $1\le i\le n$. It suffices to exclude the
latter. Suppose that $-\L_\star$ is given by~\eqref{exce}. Then
$\det(-\L_\star)=-1$, that is,
\begin{equation}\label{bi}
  \prod_{i=1}^n b_i=(-1)^n 2.
\end{equation}
Using~\eqref{redvi}, we obtain
\[
  H^*(M)=\Z[x_1,\dots,x_n]\bigl/\bigl(x_1(x_1+b_1x_2),x_2(x_2+b_2x_3),
  \ldots,x_n(b_nx_1+x_n)\bigr),
\]
where we set $x_i=v_{i+n}$ for $1\le i\le n$. By the assumption,
there is a non-zero element $x\in H^2(M,\Q)$ whose square is zero.
Write $x=\sum_{i=1}^n a_ix_i$ for some $a_i\in\Q$, then
\begin{multline*}
  0=\Bigl(\sum_{i=1}^n a_ix_i\Bigr)^2=
  \sum_{i=1}^n a_i^2x_i^2+2\sum_{i<j}a_ia_jx_ix_j\\
  =-a_1^2b_1x_1x_2-a_2^2b_2x_2x_3-\ldots-a_n^2b_nx_nx_1
+2\sum_{i<j}a_ia_jx_ix_j,
\end{multline*}
which implies that
\begin{equation}\label{aibj}
  a_1^2b_1=2a_1a_2,
  \quad a_2^2b_2=2a_2a_3, \quad\ldots,\quad a_n^2b_n=2a_na_1.
\end{equation}
Suppose $a_i\ne0$ for every $i$. Multiplying the identities above,
we obtain $\prod_{i=1}^n b_i=2^n$, which contradicts~\eqref{bi}.
Therefore $a_i=0$ for some $i$, but this together
with~\eqref{aibj} implies that $a_i=0$ for every $i$. This
contradicts the assumption that $x=\sum_{i=1}^n a_ix_i\ne0$.
Therefore, \eqref{exce} cannot occur as a reduced characteristic
matrix, and $M$ is a Bott tower.
\end{proof}

We are now ready to prove our final promised result.

\begin{theo}\label{qtcpn}
A quasitoric manifold $M$ is homeomorphic to $(\C P^1)^n$ if and
only if $H^*(M)\cong H^*((\C P^1)^n)$ as rings.
\end{theo}
\begin{proof}
Since $H^*((\C P^1)^n;\Z/2)$ is a $\BQ$-algebra of rank $n$, the
quotient polytope of $M$ is an $n$-cube by Theorem~\ref{cube}.
Then $M$ is a Bott tower by Lemma~\ref{Qcoef}. Finally, $M$ is
homeomorphic to $(\C P^1)^n$ by Theorem~\ref{trivc}.
\end{proof}

We finish by proposing the following quasitoric analogue of
Problem~\ref{cohcl}.

\begin{prob}\label{cclqt}
Does an isomorphism of rings $H^*(M_1)\cong H^*(M_2)$ imply a
homeomorphism of quasitoric manifolds $M_1$ and $M_2$?
\end{prob}

\end{document}